\documentclass[12pt]{extarticle}

\usepackage[margin=1.4in]{geometry}
\usepackage{amsthm}
\usepackage{amsmath}
\usepackage{amssymb}
\usepackage{amsfonts}
\usepackage{latexsym}
\usepackage[autostyle]{csquotes}
\usepackage{stackengine}
\usepackage{enumitem}
\usepackage{hyperref}
\usepackage{tikz-cd}
\usepackage{titlesec}
\usepackage{faktor}

\usepackage[style=alphabetic, sorting=ynt, maxbibnames=99, backend=biber]{biblatex}

\hypersetup{colorlinks,citecolor=black,filecolor=black,linkcolor=black,urlcolor=black}
\titleformat{\section}[hang]{\normalfont\large\bfseries\centering}{§\,\thesection}{1em}{}
\titleformat{\subsection}[hang]{\normalfont\bfseries}{§\,\thesubsection}{1em}{}
\newcommand\blfootnote[1]{%
	\begingroup
	\renewcommand\thefootnote{}\footnote{#1}%
	\addtocounter{footnote}{-1}%
	\endgroup
}

\newcommand{\ptsabove}{1em}
\newcommand{\ptsbelow}{1em}
\newcommand{\ptsindentamount}{0em}
\newcommand{\theoremmarin}{0em}
\newcommand{\emsbelowtheoremhead}{0.5em}

\newtheoremstyle{default theorem style}{\ptsabove}{\ptsbelow}{\addtolength{\leftskip}{\theoremmarin}\itshape}{\ptsindentamount}{\bf}{.}{\emsbelowtheoremhead}{(\thmnumber{#2}) \thmname{#1}}
\newtheoremstyle{theorem style for definitions}{\ptsabove}{\ptsbelow}{}{\ptsindentamount}{\bf}{.}{\emsbelowtheoremhead}{(\thmnumber{#2}) \thmname{#1}}

\theoremstyle{default theorem style}
\newtheorem{theorem}{Theorem}[section]
\newtheorem{lemma}[theorem]{Lemma}

\newtheorem{prop}[theorem]{Proposition}
\newtheorem{cor}{Corollary}[theorem]

\theoremstyle{theorem style for definitions}

\newtheorem{examples}[theorem]{Examples}

\setlist[enumerate,1]{label={(\roman*)}}

\newcommand{\spec}{\mathrm{Spec}}
\newcommand{\quo}{\textrm{Frac}}

\newcommand{\Q}{{\mathbb Q}}
\newcommand{\m}{\mathfrak{m}}
\newcommand{\p}{\mathfrak{p}}
\newcommand{\seq}{\subseteq}
\renewcommand{\>}{\rangle}

\renewcommand{\tilde}{\widetilde}

\title{Straight Domains are Locally Divided}
\author{Spencer Secord\footnote{This paper was written under the supervision of David McKinnon, as part of the Pure Mathematics Masters program at the University of Waterloo.}}
\date{\vspace{-2ex}}

\newcommand{\MainCorText}{A domain is straight if and only if it is locally divided.}
\newcommand{\MainTheoremText}{A prime ideal $\p$ of a domain $A$ is straight, if and only if, $\p$ is a locally divided prime ideal of $A$.}

\setlist[enumerate,1]{label={(\roman*)}}
\begin{document}
\maketitle
\vspace{-1em}
\begin{abstract}
	We present a proof that all straight domains are locally divided\textemdash thereby answering two open problems presented in \cite{open_questions}, which originally appeared in \cite{dobbs2009straight} and \cite{dobbs2009straight_II}. 
	In fact, we are able to prove a stronger result: a prime ideal of a domain is straight if and only if it is locally divided. 
\end{abstract}
\blfootnote{\textit{Date:} \today}
\blfootnote{\textit{Email:} spencer.e.secord+math@gmail.com}
\vspace{-2em}
\section{Introduction}
All rings are assumed to be unital and commutative. 
An overring $B$ of a domain $A$ is any ring extension $A\seq B$ such that $B$ is contained in the field of fractions of $A$. 
The connection between ideal theory and overring theory characterises much of modern commutative domain theory\textemdash motivated by the numerous characterisations of Dedekind Domains, Valuation Domains, and Pr\"ufer Domains \cite{Gilmer_multiplicative_ideal_theory}. 
Moreover, the definitions of most object we study in this paper (locally divided domains, CPI-extensions, and straight domains) are all generalizations of particular characterisations of Pr\"ufer domains.

A prime ideal $\p$ of a domain $A$ is a \textit{divided} prime ideal if $\p = \p A_\p$.
If all the prime ideals of $A$ are divided, then we say $A$ is a \textit{divided domain}.
Similarly, a prime ideal $\p$ is a \textit{locally divided} prime ideal of $A$ whenever $\p A_\m$ is a divided prime ideal of $A_\m$ for every maximal ideal $\m$ of $A$ that contains $\p$. 
We call $A$ \textit{locally divided} if $A_\m$ is divided for every maximal ideal $\m$ of $A$; or equivalently, if every prime ideal of $A$ is locally divided.

Divided domains were first studied in \cite{akiba1967note} as a generalization of valuation domains, 
although the naming and systematic study of divided and locally divided domains does not appear in the literature until Dobbs' papers: \cite{dobbs1976divided} and \cite{Dobbs_1981}. 
In the years since, divided and locally divided domains have been extensively studied by many other mathematicians. 
Locally divided domains are key examples of going down domains \cite[Remark 2.7 (b)]{dobbs1976divided}, and the property of being \enquote{locally divided} has been used as one motivating example in the emerging study of properties of properties in commutative ring theory (see \cite{dobbs2014new}, \cite{dobbs2009straight}).
Moreover, divided domains have inspired the study of many other interesting objects (see \cite{badawi1995domains}, \cite{badawi2002powerful}, \cite{badawi2002pseudo}, \cite{badawi2004divided}, \cite{badawi2007}).

\begin{examples}
	All maximal ideals are locally divided prime ideals, but a maximal ideal is divided if and only if the ring is local. 
	One well-studied class of locally divided domains are \textit{Pr\"ufer domains}, which are domains $A$ such that $A_\m$ is a valuation domain for every maximal ideal $\m$.
	Even though Pr\"ufer domains (and thus Dedekind domains) are locally divided, the converse is not true. 
	For instance, take $A=\Q[x,y]+\<x\>\Q[x,y]_{\<x\>}$ where $\Q[x,y]_{\<x\>}$ is the localization of $\Q[x,y]$ by the prime ideal $\<x\>$.
	This is locally divided since $A$ has Krull dimension 2, $x\Q[x,y]_{\<x\>}$ is a divided prime ideal of $A$, and every maximal ideal is locally divided. 
	To show that $A$ is not a Pr\"ufer domain, we need to show that there is a maximal ideal $\m$ such that $A_\m$ has two incomparable ideals. 
	Indeed, let $\m=yA+x\Q[x,y]_{\<x\>}$, which is the maximal ideal corresponding to the maximal ideal $\<x,y\>$ of $\Q[x,y]$. 
	Consider the ideals $\<xy^2\>$ and $\<x(x+y)\>$ of $A_\m=\Q[x,y]_{\<x,y\>}+x\Q[x,y]_{\<x\>}$. We know that $(x+y)/y^2$ and $y^2/(x+y)$ are not in $A_\m$, so $\<xy^2\>$ and $\<x(x+y)\>$ are incomparable ideals of $A_\m$. 
	Hence $A_\m$ is not a valuation domain, and therefore $A$ cannot be a Pr\"ufer domain. 
\end{examples}

We say $\p$ is a straight prime ideal of $A$ (or $A$ is \textit{straight at $\p$}) if for every overring $B$ of $A$, the ($A/\p$)-module $B/\p B$ is torsion-free.
We call $A$ a \textit{straight domain} if all its prime ideals are straight. 
The study of straight domains began with Dobbs' and Picavet's paper \cite{dobbs2009straight} which was continued in \cite{dobbs2009straight_II}. 
In these papers, the question of whether all straight domains are locally divided was first presented and studied \cite[§1]{dobbs2009straight}, which eventually led to the question of whether being a straight domain is a portable property \cite[§2]{dobbs2009straight_II}. 
Interest in this question helped motivate the systematic study of portable properties \cite{dobbs2014new}, although we will not define portable properties here due to their cumbersome characterization and lack of use in this paper. 
These two open questions on straight domains gained more interest after their appearance in \cite[Problems 13a and 13b]{open_questions}.

The main result of this paper is Theorem (\ref{main theorem}), which states:
\begingroup
\def\thetheorem{\ref{main theorem}}
\begin{theorem}
	\MainTheoremText
\end{theorem}
\addtocounter{theorem}{-1}
\endgroup
This answers both aforementioned open questions appearing in \cite{dobbs2009straight}, {\cite{dobbs2009straight_II}}, and \cite[Problems 13a and 13b]{open_questions}.
Our result is surprising since this does not hold for their generalizations; i.e. there are straight rings with zero divisors which are not locally divided rings \cite[Example 4.5]{dobbs2009straight}.
We are able to prove Theorem (\ref{main theorem}) by using CPI-extensions, which we will define shortly.	

The CPI-extension of a domain $A$ by a prime ideal $\p$ is the ring $A+\p A_\p$. 
CPI-extensions were first introduced by Boisen and Sheldon \cite{Boisen_1977}, who chose the name as an abbreviation of \enquote{complete pre-image,} which was named as such because the CPI-extension of $A+\p A_\p$ is the preimage of $A/\p$ under the canonical $A_\p\rightarrow A_\p/\p A_\p$.
A CPI-extension of a domain is always an overring, and there is a natural correspondence between prime ideals of $A$ which are comparable to $\p$ and prime ideals of the CPI-extension $A+\p A_\p$. 
CPI-extensions are used extensively in the study of divided and locally divided domains. For instance, $A$ is a divided domain if and only if every CPI-extension of $A$ is equal to $A$, and a domain is locally divided if and only if every CPI-extension is a localization \cite{Dobbs_1981}.

\subsection{Summary}
Section \ref{Section: Preliminary Lemmas} contains the prerequisite results needed to prove our main result, Theorem (\ref{main theorem}).
Section \ref{Section: the main results} focuses on proving Theorem (\ref{main theorem}) and its corollaries. 
In particular, Theorem (\ref{main theorem}) states: \enquote{\textit{\MainTheoremText}}
We then show that this answers two open questions appearing in {\cite{dobbs2009straight}}, \cite{dobbs2009straight_II}, and \cite[Problems 13a and 13b]{open_questions}.
Namely, that every straight domain is locally divided (\ref{main cor}), and that \enquote{being a straight domain} is a portable property of domains (\ref{main cor 2}).

\section{Preliminary Lemmas}\label{Section: Preliminary Lemmas}

We start by proving a slight generalization of a result appearing in Kaplanski's book \cite{Kaplanski}.

\begin{prop}\textup{\cite[Theorem 67]{Kaplanski}.}\label{Seidenberg's Lemma}
	Let $A\seq B$ be a ring extension with $A$ local and integrally closed in $B$. 
	Suppose $f(x)\in A[x]$ is a polynomial with atleast one unit coefficient. 
	If $b\in B$ such that $f(b)=0$, then either $b\in A$ or $b$ is invertible in $B$ and $b^{-1}\in A$.
	\begin{proof}
		We will follow the proof of \cite[Theorem 67]{Kaplanski}, generalizing where needed.
		Since $A$ is integrally closed in $B$, we may assume the leading coefficent of $f=\sum_{i=0}^n a_i x^i$ is not a unit. 
		We will prove this by induction on the degree of $f$.
		For the base case $a_1b-a_0=0$, we know $a_1b$ is a unit of $A$, so $b$ must be invertible in $B$. 
		This implies $b^{-1}=\frac{a_1}{a_0}\in A$, as needed.
		
		Now for the induction step.
		Notice that $(a_nb)^n+\sum_{i<n} a_i a_n^{n-i} (a_nb)^i=0$ and thus $a_nb\in A$. 
		Therefore 
		\((a_nb-a_{n-1})b^{n-1}+\sum_{i<n-1}a_i b^i=0\).
		If the unit coefficient of $f$ is not $a_{n-1}$, then we can apply the induction assumption to $(a_nb-a_{n-1})x^{n-1}+\sum_{i<n}a_i x^i$ and we are done.
		If $a_{n-1}$ is a unit, then since $a_n$ is not a unit and $A$ local, we know $a_nb-a_{n-1}$ must be a unit. 
		Since $A$ is integrally closed in $B$, we would have $b\in A$, as needed.
	\end{proof}
\end{prop}

The following lemma will be used to prove Theorem (\ref{main theorem}), which will make this lemma redundant.

\begin{lemma}\label{p overring lemma}
	Let $A$ be a local domain which is straight at the prime ideal $\p$. 
	\begin{enumerate}[label=(\roman*)]
		\item \label{p overring lemmas lemma 1} If $b\in \p A_\p$ then $b\in \p A[b]$.
		\item \label{straight overring is integral in CPI-ext lemma}$A$ is integrally closed in $A+\p A_\p$.
	\end{enumerate}
	\begin{proof}
		First lets prove (i).
		Without loss of generality, let $x\in\p$ and $y\in A\setminus \p$. 
		Then we have $y\frac{x}{y} = x \in \p \seq \p A[\frac{x}{y}]$. 
		Since $A$ is straight at $\p$ we know $\frac{x}{y}\in\p A[\frac{x}{y}]$, as needed.
		
		Next we will prove (ii).
		Suppose $b\in A+\p A_\p$ is integral over $A$.
		Without loss of generality, we may assume $b\in \p A_\p$, and thus $b=\frac{x}{y}$ for some $x\in \p$ and $y\in A\setminus\p$. 		
		We know $A[\frac{x}{y}]$ is an overring of the straight domain $A$ and $y\frac{x}{y} = x\in\p\seq \p A[\frac{x}{y}]$, so $\frac{x}{y}$ must be in $\p A[\frac{x}{y}]$. 
		Since $\p A[\frac{x}{y}]$ is multiplicatively closed, we have $\frac{x^i}{y^i}\in \p A[\frac{x}{y}]$ for every $i$. 
		Therefore $\p A[\frac{x}{y}] + \frac{x}{y} A[\frac{x}{y}] = \p A[\frac{x}{y}]$. 
		Since $\p A[\frac{x}{y}]=\p+\p\frac{x}{y} A[\frac{x}{y}]$ it follows that $\p + \frac{x}{y}A[\frac{x}{y}] = \p +\p \frac{x}{y}A[\frac{x}{y}]$.
		
		Let $M = \faktor{\left(\p + \frac{x}{y}A[\frac{x}{y}] \right)}{\p}$. By the above argument, we have
		\[\p M = \p \left(\faktor{\p + \frac{x}{y}A\left[\frac{x}{y}\right]\; }{\p}\right) = \faktor{ \p +\p \frac{x}{y} A\left[\frac{x}{y}\right]\; }{\p} = \faktor{\p + \frac{x}{y} A\left[\frac{x}{y}\right]\; }{\p} = M\]
		Now let $J$ be the Jacobson radical of $A$. Since $A$ is local, $J$ must be the maximal ideal $\m$ of $A$, and thus must contain $\p$. Hence $JM = M$. 
		Recall that $\frac{x}{y}$ is integral over $A$, which implies $A[\frac{x}{y}]$ must be a finitely generated $A$-module, and therefore $M$ must be a finitely generated $A$-module.		
		Since $M$ is a finitely generated $A$-module with $J M = M$, we can apply Nakayama's Lemma to get $M=0$, and thus $\frac{x}{y}A[\frac{x}{y}]\seq \p\seq A$ which implies $\frac{x}{y}$ is in $A$. Therefore, $A$ must be integrally closed in $A+\p A_\p$.
	\end{proof}
\end{lemma}

\section{The Main Results}\label{Section: the main results}
We can now prove the main result of this paper:

\begin{theorem}\label{main theorem}
	\MainTheoremText
	\begin{proof}
		Note that the property \enquote{\itshape $\p A_\p = \p A_\m$ for every maximal ideal $\m$ that contains $\p$} is a local property over domains. 
		Let us show that $A$ is straight at $\p$ if and only if $A_\m$ is straight at $\p_\m$ for every maximal ideal $\m$ \cite[Proposition 3.5]{dobbs2009straight}.
		For an overring $B$ of $A$, the property \enquote{$B/\p B$ is a torsion-free $(A/\p)$-module} is a local property with respect to $\spec(A/\p)$. Since the localization of $B/\p B$ by a maximal ideal of $A$ which does not contain $\p$ must be trivial, we know that \enquote{$B/\p B$ is a torsion-free $(A/\p)$-module} is a local property with respect to the maximal ideals of $A$. Therefore, \enquote{$A$ is straight at $\p$} must be a local property.
		It is therefore sufficient to prove this theorem when $A$ is a local domain.
		
		Assume $A$ is local and let $b\in \p A_\p$, we will show $b\in A$. 
		By \ref{p overring lemmas lemma 1} of Lemma (\ref{p overring lemma}), we know $b\in \p A[b]$, and thus $b=\sum_i a_i b^i$. 
		Define the polynomial $f(x)= a_0+(a_1+1)x+\sum_{i>1}a_i x^i$, and notice that $f(b)=0$. 
		Clearly the $x$ coefficient is a unit since $A$ is local and $a_1$ is not a unit.
		By Proposition (\ref{Seidenberg's Lemma}), we know either $b\in A$ or $b^{-1}\in A$.
		Yet $b$ is a non-unit in $A+\p A_\p$, so $b^{-1}\not\in A$. 
		This implies $b\in A$, and therefore $A=A+\p A_\p$.		
	\end{proof}
\end{theorem}

The following corollary answers the question posed in \cite[§1]{dobbs2009straight}, which gained further interest after appearing in \cite[Problems 13b]{open_questions}.

\begin{cor}\label{main cor}
	\MainCorText
	\begin{proof}
		Clearly every locally divided prime ideal is straight, and thus the statement follows trivially from Theorem (\ref{main theorem}).
	\end{proof}
\end{cor}
The next corollary answers the question posed in \cite[§2]{dobbs2009straight_II}.
This question also appeared in \cite[Problems 13a]{open_questions}, under an equivalent but different phrasing.
Although the following corollary is about portable properties, we have omitted the definition of a portable property since it is cumbersome and will not be used here. See \cite{dobbs2009pullback} for a detailed description.
\begin{cor}\label{main cor 2}
	\enquote{\itshape Being a straight domain} is a portable property.
	\begin{proof}
		\cite[Lemma 2.1]{dobbs2009pullback} proves \enquote{\itshape being a locally divided domain} is a portable property over domains. By Corollary (\ref{main cor}) we are done.
	\end{proof}
\end{cor}
Finally, we can use Theorem (\ref{main theorem}) to characterize locally divided prime ideals of domains.
\begin{cor}\label{straight at all domain extensions}
	Let $\p$ be a prime ideal of a domain $A$. Then the following are equivalent:
	\begin{enumerate}[label=(\roman*)]
		\item $\p$ is a locally divided prime ideal of $A$; i.e. $\p A_\m = \p A_\p$ for every prime ideal $\m$ containing $\p$.
		\item $\p$ is a straight prime ideal of $A$; i.e. for every overring $B$ of $A$ the $(A/\p)$-module $B/\p B$ is torsion-free.
		\item For every domain extension $(A, B)$, the $(A/\p)$-module $B/\p B$ is torsion-free.
	\end{enumerate}
\begin{proof}
	By Theorem (\ref{main theorem}) we know that (i) is equivalent to (ii). Clearly (iii) implies (ii), so we need only show that (i) and (ii) imply (iii).
	
	Let $A\seq B$ be an arbitrary domain extension and let $K=\quo(B)$ be the field of fractions of $B$. Now consider the following subrings of the field of rational functions $K(x)$ given by:
	\[\tilde A = A+ x K[x]_{\<x\>},\quad \tilde B = B+ x K[x]_{\<x\>}\]
	where $K[x]_{\<x\>}$ is the localization of $K[x]$ by the prime ideal $\<x\>$. Note that $(\tilde A,\tilde B)$ is an overring extension with field of fractions $K(x)$. 
	
	Let $\tilde\p=\p+x K[x]_{\<x\>}$. Clearly, $\tilde A/\tilde \p$ is naturally isomorphic to $A/\p$, and through this isomorphism we know $B/\p B$ is naturally isomorphic to $\tilde B/\tilde \p \tilde B$ as $(A/\p)$-modules. Thus it suffices to show that $\tilde\p$ is divided, yet $x K[x]_{\<X\>}$ is divided and $\p\cong\tilde\p/x K[x]_{\<x\>}$ is a locally divided prime ideal, and therefore $\tilde \p$ must be locally divided. Since (i) is equivalent to (ii), we know that $B/\p B\cong\tilde B/\tilde \p \tilde B$ must be torsion-free as an $(A/\p)$-module.
\end{proof}
\end{cor}

\section{Acknowledgements}\label{Section: Acknowledgements}
I would like to thank my supervisor David McKinnon for his mentorship, insightful advice, and contagious curiosity. 
I am thankful towards Jason Bell for carefully reading this paper and for providing numerous insightful comments, suggestions, and examples. 
I am also thankful towards Thomas Bray for his frequent help and many insightful discussions on mathematical topics, and to Dan Labach for proofreading this paper.

I am grateful to Tiberiu Dumitrescu for his help in creating the final form of this paper: Shortly after submitting my preprint to Arxiv, Dumitrescu notified me of an error in the first preprint of this paper, which related to an incorrect application of a local-global argument in a lemma. By instead using Proposition \ref{Seidenberg's Lemma}, Dumitrescu was able to remedy this problem and generously sent me his work, along with permission to use it in this document.

\newpage
\printbibliography
\end{document}